\input amstex.tex
\documentstyle{amsppt}
\input amssym.def
\newsymbol\geqq 133D
\TagsOnRight
\newsymbol\blacksquare 1004

\def\bn{{\Bbb N}}
\def\br{{\Bbb R}}
\def\bc{{\Bbb C}}

\def\bz{{\Bbb Z}}

\def\X{{\bar X}}

\def\cB{{\Cal B}}

\def\cL{{\Cal L}}

\def\cH{{\Cal H}}
\def\cO{{\Cal O}}

\def\sig{{\sigma}}
\def\Om{{\Omega}}
\def\lam{{\lambda}}
\def\Lam{{\Lambda}}

\def\hm{{\hat\mu}}

\def\Mat{\operatorname{Mat}}
\def\tr{\text{tr}}
\def\lat{\operatorname{lat}}
\def\andd{\qquad\text{and}\qquad}
\magnification=1100
\baselineskip=2\normalbaselineskip
\hfuzz=5pt 

\topmatter
\title Harmonic Analysis of Fractal Measures \endtitle

\rightheadtext{Harmonic Analysis of Fractal Measures}
\author Palle E.T. Jorgensen and Steen Pedersen \endauthor
\address Department of Mathematics, University of Iowa, Iowa City, IA  52242, USA 
\endaddress
\address Department of Mathematics, Wright State University, Dayton, OH  45435, USA 
\endaddress
\thanks Research supported by the NSF. \endthanks

\abstract
We consider affine systems in $\br^n$ constructed from a given integral invertible and 
expansive matrix $R$, and a finite set $B$ of translates, $\sig_b x:=R^{-1}x+b$; the 
corresponding measure $\mu$ on $\br^n$ is a probability measure and fixed by the 
selfsimilarity $\mu=|B|^{-1}\sum_{b \in B} \mu \circ \sig_b^{-1}$. There are two {\it a 
priori\/} candidates for an associated orthogonal harmonic analysis : (i) the existence of 
some subset $\Lam$ in $\br^n$ such that the exponentials $\{e^{i\lam\cdot x}\}_{\lam \in 
\Lam}$ form an {\it orthogonal basis\/} for $L^2(\mu)$; and (ii) the existence of a certain 
{\it dual pair of representations\/} of the $C^*$-algebra $\cO_N$ where $N$ is the 
cardinality of the set $B$. (For each $N$, the $C^*$-algebra $\cO_N$ is known to be 
simple; it is also called the Cuntz-algebra.) We show that, in the ``typical'' fractal case, the 
naive version (i) must be rejected; typically the orthogonal exponentials in $L^2(\mu)$ fail 
to span a dense subspace. Instead we show that the $C^*${\it -algebraic version of an 
orthogonal harmonic analysis\/}, viz., (ii), is a natural substitute. It turns out that this 
version is still based on exponentials $e^{i\lam\cdot x}$, but in a more indirect way. (See 
details in Section 5 below.) Our main result concerns the intrinsic geometric features of 
affine systems, based on $R$ and $B$, such that $\mu$ has the $C^*$-algebra property 
(ii). Specifically, we show that $\mu$ has an orthogonal harmonic analysis (in the sense 
(ii)) if the system $(R,B)$ satisfies some specific symmetry conditions (which are 
geometric in nature). Our conditions for (ii) are stated in terms of two pieces of data: (a) a 
{\it unitary generalized Hadamard-matrix\/}, and (b) a certain {\it system of lattices\/} 
which must exist and, at the same time, be compatible with the Hadamard-matrix. A 
partial converse to this result is also given. Several examples are calculated, and a new 
maximality condition for exponentials is identified.
\endabstract
\endtopmatter

\document
\head 1. Introduction \endhead
\openup 1.0\baselineskip
The present paper continues work by the coauthors in \cite{JP3--6}, and it also provides 
detailed proofs of results announced in \cite{JP4}. In addition we have new results going 
beyond those of the announcement \cite{JP4}. We consider a new class of selfsimilarity 
fractals $\X$, each $\X$ with associated fractal selfsimilar measure $\mu$, such that 
$L^2(\mu)$ has an orthogonal harmonic analysis in the sense of $C^*$-algebras (see (ii) 
below). This possibility is characterized with geometric axioms on the pair $(\X, \mu)$. It 
is known since \cite{St3--4} that $\mu$ is typically singular (in the fractal case), and that 
in general only an {\it asymptotic\/} Plancherel type formula can be expected in the sense 
of \cite{Bes}. Our present approach is based instead on $C^*$-algebra theory. In 
particular, we use the $C^*$-algebras $\cO _N$ of Cuntz \cite{Cu}, and we give the 
orthogonal decompositions in terms of a {\it dual pair of representations\/} of $\cO _N$ 
where $N$ denotes the number of translations in the affine system which determines 
$\mu$. 

For an orthogonal harmonic analysis, the following {\it three possibilities\/} appear {\it a 
priori\/} as natural candidates: 
\item{(i)} the existence of a subset $\Lam$ in $\br^n$ such that the exponentials $e_\lam 
(x):=e^{i\lam\cdot x}$ (indexed by $\lam \in \Lam$) form an {\it orthogonal basis\/} in 
$L^2(\mu)$;
\item{(ii)} the existence of a {\it dual system of representations\/} of some $C^*$-algebra 
$\cO_N$ say, ($N=$ the cardinality of $B$), such that one representation is acting affinely 
in $x$-space, and the other (dually) in frequency-space (where the frequency variable is 
represented by $\lam$ in the above exponentials $e^{i\lam\cdot x}$); and finally
\item{(iii)} one might base the harmonic analysis on an orthogonal basis of polynomials in 
$n$ variables obtained from the monomials $x^{\alpha}:= x_1^{\alpha_1} 
x_2^{\alpha_2}\cdots x_n^{\alpha_n}$ (where $\alpha=(\alpha_1, \ldots,\alpha_n$) is a 
multi-index, $\alpha_i =0,1,2,\ldots$, $1\leq i \leq n$), by the familiar Gram-Schmidt 
algorithm.

But it is immediate that both of the possibilities (i) and (iii) lack symmetry in the variables 
$x$ and $\lam$. Moreover, it turns out that (i) must be ruled out also for a more serious 
reason. We show in Sections 6--7 below that, for the ``typical'' fractal measures $\mu$, 
none of the orthogonal sets $\{e_\lam\}$ in $L^2(\mu)$ will in fact span a dense 
subspace. Specifically, there is a canonical maximally orthogonal $\{e_\lam\}$ system 
such that a finite set of ``translates'' (details in Section 5) of it does give us a dense 
subspace. It is this extra operation (i.e., ``spreading out'' the orthogonal exponentials) 
which leads to our dual pair of representations of the algebra $\cO_N$. 

It also turns out that case (ii) is a natural extension of our orthogonality condition, studied 
earlier in \cite{JP2} for $L^2(\Om)$, now with $\Om$ some subset in $\br^n$ with finite 
positive Lebesgue measure, and $L^2(\Om)$ considered as a Hilbert space with the 
restricted Lebesgue measure. For the case, when $\Om$ is {\it further assumed\/} open 
and connected, we showed, in \cite{JP2} and \cite {Pe}, that (i) holds (i.e., there is a set 
$\Lam$ such that $\{e_\lam \}_{\lam \in \Lam}$ forms an orthogonal basis in $L^2 
(\Om)$) iff the corresponding {\it symmetric\/} operators $\{\sqrt{-
1}\frac{\partial}{\partial x_j} \}^n _{j=1}$, defined on $C^\infty _c (\Om)$, have 
commuting {\it selfadjoint\/} operator extensions acting in $L^2(\Om)$. It is well known 
that extension theory for symmetric operators is given by von Neumann's deficiency 
spaces. But, even when individual selfadjoint extensions exist for commuting symmetric 
operators, such extensions are typically non-commuting. Hence, we expect that, also for 
our $L^2(\mu)$ analysis, there will be {\it distinct\/} symmetry conditions and 
selfadjointness conditions.

For our present case, the pair $(R,B)$ is specified as above, the affine maps are given by 
$\sig_b x=R^{-1}x+b$, and indexed by points $b$ in the finite set $B$. We get the 
measure $\mu$, and the Hilbert space $L^2(\mu)$, by a general limit construction which 
we show must start with some $L^2(\Om)$ example as discussed. But, for $L^2(\mu)$, 
we show that the analogous {\it symmetry\/} condition is related to a certain lattice 
configuration in $\br^n$ (see Lemma 4.1 below), whereas the analogous {\it 
selfadjointness\/} now corresponds to a spectral pairing between $B$ and a second subset 
$L$ in $\br^n$, of same cardinality, such that the $N$ by $N$ matrix $\{e^{ib 
\cdot\ell}\}$, (for $b\in B$, $\ell \in L$), forms a so-called {\it unitary generalized 
Hadamard matrix\/}, see \cite{SY}. Then this matrix, together with the lattice 
configuration leads to a dual pair of representations, as sketched above and worked out in 
detail below. The two representations will act naturally on $L^2(\mu)$ and provide a non-
commutative harmonic analysis with a completely new interpretation of the classical time-
frequency duality (see e.g., \cite{HR}), of multivariable Fourier series.

When our ``symmetry'' condition is satisfied, we get a dual pair of self-similar measures, 
$\mu_B$ and $\mu_L$, and this pair is used in the proof of our structure theorem. Many 
examples are given illustrating when the ``symmetry'' holds and when it doesn't. A 
connection is made to classical spectral duality, see e.g., \cite{JP1--3}.

\head 2. Basic Assumptions \endhead
We consider affine operations in $\br^n$ where $n$ is fixed; the case $n=1$ is also 
included, and the results are non-trivial and interesting also then. A {\it system\/} $s$ in 
$\br^n$ will consist of a quadruple $(R,B,L,K)$ where $R\in GL_n(\br)$, $B$ and $L$ 
are {\it finite subsets\/} in $\br^n$, and both of them are assumed to contain the origin 
$O$ in $\br^n$; finally $K$ is a {\it lattice\/} in $\br^n$, i.e., a free additive group with 
$n$ generators. It will be convenient occasionally to identify a fixed {\it lattice\/} with a 
{\it matrix\/} whose columns are then taken to be a set of generators for the lattice in 
question. It is known that generators will always form a linear basis for the vector space 
(see e.g., \cite{CS}); and it follows that the matrix is then in $GL_n(\br)$.

With the assumptions (to be specified), it turns out that we may apply Hutchinson's 
theorem \cite{Hu} to the affine system $\{\sig_b \}_{b \in B}$ given by 
$$
\sig_b x:= R^{-1} x + b,\qquad x \in\br^n \tag 2.1
$$

There is a unique probability measure $\mu$ on $\br^n$ satisfying
$$
\mu = |B|^{-1} \sum_{b \in B} \mu \circ \sig_b^{-1},\tag 2.2
$$
which amounts to the condition
$$
\int f\,d\mu=|B|^{-1} \sum_b \int f \circ \sig_b\, d\mu \tag 2.3
$$
for all $\mu$-integrable functions $f$ on $\br^n$. For the matrix $R$, we assume that 
some positive integral power of it has all eigenvalues in $\{\lam \in \bc : |\lam|>1\}$, and 
we refer to this as the {\it expansive\/} property for $R$. (It is actually equivalent to the 
same condition for $R$ itself.) The use of \cite{Hu} requires the so called {\it open-set-
condition\/} which turns out to hold when our system $s$ has a {\it symmetry\/} property 
which we proceed to describe. We then also have the following compact subset $\X$, 
defined as the closure (in $\br^n$) of the set of vectors $x$ with representation
$$
x=\sum_{i=0}^\infty R^{-i}b_i,\qquad b_i \in B.\tag2.4
$$
If $|B|<|\det R|$, where $|B|$ denotes the cardinality of $B$, then the 
{\it fractal dimension\/} of $\X$ will be less than the vector space dimension $n$ of the 
ambient $\br^n$. (See e.g., \cite{Ke} for details on this point.) In general, the measure 
$\mu$ is supported by $\X$, and we may identify $L^2(\mu)$ with $L^2(\X,\mu)$ as a 
{\it Hilbert space\/}. We will refer to $\X$ as the ``{\it fractal\/}'' even in the cases when 
its dimension may in fact be integral, and the ``fractal'' representation will be understood to 
be (2.4). Occasionally, we will write $\X(B)$ to stress the digit-set $B$.

\head 3. Generalized Hadamard Matrices \endhead
The two sets $B$ and $L$ from the system came up in our previous work (see \cite{JP2--
4} and \cite{JP6}) on {\it multivariable spectral theory\/}. The condition we wish to 
impose on two sets $B,L$ amounts to demanding that the corresponding exponential 
matrix
$$
\left( e^{i2\pi b\cdot \ell} \right) \tag 3.1
$$
is {\it generalized Hadamard\/}, see \cite{SY}. The term $b \cdot \ell$ refers to the usual 
dot-product in $\br^n$. It will be convenient to abbreviate the matrix entries as, $\langle 
b,\ell \rangle :=e^{i2\pi b\cdot \ell}$. Since $0 \in B$ and $0 \in L$ by assumption, one 
column, {\it and\/} one row, in the matrix $(\langle b,\ell \rangle )_{BL}$ consists of a 
string of ones. Let the matrix be denoted by $U$: We say that it is {\it generalized 
Hadamard\/} if the two sets $B$ and $L$ have the same cardinality, $N$ say, and if
$$
U^* U=NI_N. \tag 3.2
$$
It follows from this that then also $UU^*=NI_N$. (This is just saying, of course, that the 
complex $N$ by $N$ matrix, $N^{-1/2}U$ is unitary in the usual sense.)

We noted in \cite{JP6} that the harmonic analysis of type (ii) is based on this kind of 
Hadamard matrices. (The matrices also have an independent life in combinatorics.) It turns 
out that the matrices are known for $N$ up to $N=4$. We will show, in Section 7 below, 
that this then leads to a classification of the simplest affine fractals (as specified) such that 
the analysis (ii) exists. We say that two matrices $U$ of the form (3.2) are {\it 
equivalent\/}, if $N$ is the same for the two matrices, and if one arises from the other by 
multiplication on the left, or right, with a permutation matrix, or with a unitary diagonal 
matrix. We now list below (without details) the {\it inequivalent\/} cases of type (3.2) for 
$N \leq 4$. (For higher $N$, such a classification is {\it not\/} known.) After our present 
preprint was circulated, we learned that the $N\leq 4$ classification had also been found 
independently, see references \cite{Cr} and \cite{Wer}. The purpose of our examples in 
Section 7 is to show how the equivalence classes of (3.2) lead to distinct examples of 
fractal measures $\mu$, and how the different $U$-matrices lead to different dual pairs of 
representations. 

We will postpone to a later paper a rigorous classification of the different systems 
$(R,B)$, and of the corresponding type (ii) harmonic analysis of $L^2(\mu)$. But we feel 
that the $N\leq 4$ examples are sufficiently interesting in their own right. They also serve 
to illustrate the technical points in our (present) two main theorems.

Notice the $2\pi$ factor in the exponential (3.1) above. It is put in for technical 
convenience only.

\remark{Remark \rom{3.1}}If we pick the string of ones as first row and first column, 
then the possibilities for $U$ when $N=2$ are 
$$
\pmatrix1&1\\1&-1\endpmatrix ; \tag3.3
$$
for $N=3$,
$$
\pmatrix 1&1&1\\1&\zeta&\bar \zeta\\1&\bar \zeta&\zeta\endpmatrix \tag 3.4
$$
where $\zeta$ is a primitive 3rd root of 1; and for $N=4$,
$$
\pmatrix1&1&1&1\\
1&1&-1&-1\\
1&-1&u&-u\\
1&-1&-u&u \endpmatrix \tag 3.5
$$
where $|u|=1$, up to {\it equivalence\/} for generalized Hadamard matrices, see e.g. 
\cite{SY}.

\endremark

\head 4. Selfadjoint Systems \endhead
Corresponding to the affine mappings (2.1) for a given system $s=(R,B,L,K)$ we have
$$
\tau_{\ell}(t):=R^*t + \ell,\qquad t\in \br^n \tag 4.1
$$
and the inverses
$$
\tau_{\ell}^{-1} (t):=R^{*^{-1}}(t- \ell) \tag 4.2
$$
where the translations for (4.2) are given by the vectors
$$
b':=-R^{*^{-1}}(\ell),\qquad\text{as $\ell$ varies over $L$}. \tag 4.3
$$
Here $R$ is an $n$ by $n$ matrix as specified above, $B$ and $L$ are finite subsets in 
$\br^n$ both containing $O$, and $K$ is a rank~$n$ lattice. The invariance $R(K) \subset 
K$ will be assumed, and we summarize this by the notation $K \in \lat (R)$.

We introduce the {\it dual system\/} $s^\circ$ defined by $s^\circ =(R^*,B',L',K^\circ)$ 
where $K^\circ$ is {\it the dual latttice,}
$$
B':=-R^{*^{-1}}(L)\tag 4.4
$$
and
$$L':=-R(B) \tag 4.5
$$

The system $s$ is said to be {\it symmetric\/} if
$$
R(B) \subset K, \tag 4.6
$$
and if $K \in \lat (R)$; and it is said to be {\it selfadjoint\/} if both $s$ and $s^\circ$ are 
symmetric. (Also notice that, in general, we have $s^{\circ\circ}=s$ when $s$ is an 
arbitrary system.)

(The definitions are analogous to familiar ones for closed operators $S$ with dense 
domain in Hilbert space, see e.g., \cite{Fu}: the operator $S$  is said to be {\it 
symmetric\/} if $S\subset S^*$, where $S^*$ denotes the adjoint, and the inclusion refers 
to inclusion of graphs. It follows that $S$ is {\it selfadjoint\/}, i.e., $S=S^*$, iff both $S$ 
and $S^*$ are symmetric.)

We shall need the fact that $B$ embeds into the of coset space $R^{-1}(K)/K$ when {\it 
additional\/} orthogonality is assumed: 

\proclaim{Lemma 4.1}Consider a system $s=(R,B,L,K)$ in $\br^n$ with the matrix $R$, 
the two finite subsets $B$ and $L$ in $\br^n$, and a lattice $K$ as described above. 
Assume that $L\subset K^\circ$ and that the two conditions \rom{(3.2)}, i.e., that 
Hadamard property, and \rom{(4.6)} hold, where $K^\circ$ is the dual lattice in $\br^n$. 
Then it follows that different points in $B$ represent distinct elements in the finite group 
$R^{-1}(K)/K$.\endproclaim

\demo{Proof} Suppose $b\neq b'$ in $B$. Then $\sum_{\ell \in L} \langle \ell, b-b' \rangle 
= 0 $ using (3.2). But, for all $k \in K$, we also have $\sum_{\ell \in L} \langle \ell, k 
\rangle = |L|$, and it follows that $b-b' \notin K$; i.e., the $R^{-1}(K)/K$ cosets are 
distinct.\qed \enddemo

We shall assume in the following that our given system is of {\it Hadamard type\/}, i.e., 
that $|B|=|L|$ and that the matrix (3.1) formed from $(B,L)$ is generalized Hadamard, see 
(3.2) above.

The following lemma is also simple but useful.

\proclaim{Lemma 4.2} A given system $s=(R,B,L,K)$ in $\br^n$ is selfadjoint if and only 
if the following three conditions hold:\roster
\item"{(i)}"$K\in \lat(R)$,
\item"{(ii)}"$R(B)\subset K$, and
\item"{(iii)}"$L\subset K^\circ$.\endroster
\endproclaim

\demo{Proof} A calculation shows that $K \in \lat (R)$ holds iff $K^\circ \in \lat(R^*)$. 
For the system $s^\circ$ to be symmetric, we need $R^*(B') \subset K^\circ$, and that is 
equivalent to (iii) by virtue of formula (4.4). So both $s$ and $s^\circ$ are symmetric 
precisely when (i)--(iii) hold. \qed
\enddemo

\remark{Remark \rom{4.3}} (Classical Systems) In \cite{JP2}, we considered the 
following spectral problem for measurable subsets $\Om \subset \br^n$ of finite positive 
Lebesgue measure, i.e., $0<m(\Om)<\infty$ where $m=m_n$ denotes the $\br^n$-
Lebesgue-measure: Let $\Om$ be given, {\it when\/} is there a subset $\Lam\subset 
\br^n$ s.t. the exponentials
$$
e_{\lam}(x) = \langle \lam, x \rangle =e^{i2\pi \lam \cdot x} , \tag 4.7
$$
indexed by $\lam \in \Lam$, form an orthonormal basis in $L^2(\Om)$ with inner product
$$
m(\Om)^{-1} \int_{\Om} \overline{f(x)}g(x)\,dx\,? \tag 4.8
$$

The problem (in its classical form) goes back to \cite{Fu}, and it is motivated by a 
corresponding one for commuting vector fields on manifolds with boundary, see also 
\cite{Jo1--2}, \cite{Pe}, and \cite{JP2}.

We showed that the general problem may be ``reduced'' (by elimination of ``trivial'' 
systems) to a special case when the pair $(\Om,\Lam)$ is such that the polar
$$
\Lam^\circ = \{t \in \br^n: \langle t,\lam \rangle=1,\qquad \forall \lam \in \Lam\} \tag 4.9
$$
is a lattice in $\br^n$, say $K:=\Lam^\circ$, and the natural torus mapping $\br^n 
\rightarrow \br^n /K$ is then 1-1 on $\Om$.

In this case, there is a system $s=(R,B,L,K)$ which is self-adjoint and of Hadamard type. 
Moreover the set $\Lam$ (called the spectrum) may be taken as
$$
\Lam = L+R^*K^\circ. \tag 4.10
$$

Pairs $(\Om,\Lam)$ with the basis-property are called {\it spectral pairs\/}; the ``reduced'' 
ones where $\Lam$ may be brought into the form (4.10) (with $L \neq \{0\}$, i.e., 
$|L|>1$) are called {\it simple factors\/}. We showed in \cite{JP6} that more general ones 
may be built up from the simple factors.\endremark

The following easy fact will be used below: Let $K_1$ and $K_2$ be lattices, and let 
$\tilde K_1$ and $\tilde K_2$ be corresponding matrices. Then we have the lattice 
inclusion $K_1 \subset K_2$ if and only if the matrices factor: $\tilde K_1 = \tilde K_2 
M$ with $M \in \Mat_n(\bz)$ where $\Mat_n(\bz)$ denotes {\it the ring of integral $n$ by 
$n$ matrices\/}, i.e., $M=\left( m_{ij}\right)^n_{i,j=1}$ with $m_{ij}\in \bz$. This 
observation allows us to take advantage of the Noetherian property of the ring 
$\Mat_n(\bz)$. A {\it minimal\/} choice for $K$ subject to conditions is then always well 
defined.

For a given lattice $K$, the {\it dual lattice\/} is denoted $K^\circ$ and given by
$$
K^\circ :=\{s \in \br^n: s\cdot k \in \bz,\qquad \forall k \in K\} 
$$
If $\tilde K$ is a matrix for $K$, then the inverse transpose, i.e., $(\tilde K^{\tr})^{-1}$ 
will be a matrix for $K^\circ$.

When $R \in GL_n (\br)$ is given, we denote by $\lat(R)$ {\it the set of all lattices $K$ in 
$\br^n$ such that\/} $R(K)\subset K$. For the matrices, that reads $\tilde K^{-1}R \tilde 
K \in \Mat_n(\bz)$. This fact will be used in the paper; it implies for example that $|\det 
R|$ is the index of $K$ in $R^{-1}(K)$. It is known (see e.g., \cite{CS} or \cite {JP6}) 
that, if $\lat(R)\neq \emptyset$, then $\det R \in \bz$. (Remark: If $R$ is {\it not\/} in 
diagonal form, i.e., $\pmatrix  r&0&\hdots&0\\ 0&r&\hdots&0\\
\vdots&\vdots&\ddots&\vdots\\
0&0&\hdots&r \endpmatrix = rI_n$ for $r \in \bz$, then there are lattices $K$ not in 
$\lat(R)$.)

The standing assumption which is placed on $R$ is referred to as the {\it expansive\/} 
property: We assume that, for some $p \in \bn$, all the eigenvalues $\{\lam_j\}$ of $R^p$ 
satisfy $|\lam_j|>1$. Recall, $R$ has real entries, but the eigenvalues may be complex. For 
emphasis, we will denote the transpose of $R$ by $R^*$, even though it is the same as 
$R^{\tr}$. (Note that the assumption on the eigenvalues of $R^p$ for some positive 
power $p$ is equivalent to the same condition on $R$ itself, i.e., to the condition for 
$p=1$.)

\head 5. Iteration Systems \endhead
In this paper, we shall study fractals (in the sense of (2.4) above) with a high degree of 
symmetry; and show that these fractals are precisely those which may be built from 
systems $s=(R,B,L,K)$ which are {\it selfadjoint\/}, of {\it Hadamard-type\/}, and where 
the lattice $K$ is chosen as {\it minimal\/} relative to the three conditions (i)--(iii) in 
Lemma 4.1. In describing our limit systems (typically fractals), we show again that the 
Hadamard condition (3.2) is the central one.

Motivated by (4.10), we form the set $\cL (L)$ consisting of all (finite) sums
$$
\ell_0 +R^* \ell_1 + R^{*^2}\ell_2 + \cdots +R^{*^m}\ell_m \tag 5.1
$$
when $m$ varies over $\{0,1,2,\ldots\}$ and $\ell_i \in L$. Using (4.1), also notice that 
$\cL(L)$ is made from iterations
$$
\tau_{\ell_0}(\tau_{\ell_1}(\cdots(\tau_{\ell_m}(0))\cdots)). \tag 5.2
$$
The set $\Lam$ in (4.10) is $\bigcup \{\tau_\ell (K^\circ):\ell \in L\}$. We shall also need 
the corresponding iterations,
$$
\bigcup_m \{\tau_{\ell_0}\circ \cdots \circ \tau_{\ell_m}(K^\circ):\ell_i \in L\}.
\tag 5.3
$$
For a given string $(\ell_0,\ldots,\ell_m)$, the set in (5.3) will be denoted $K^\circ 
(\ell_0,\ldots,\ell_m)$.

\definition{Definition 5.1} We say that $K^\circ$ formed from a given system 
$s=(R,B,L,K)$ is {\it total\/} if the functions $\{e_s:s \in K^{\circ}\}\subset L^2(\mu)$ 
span a subspace which is dense in the Hilbert space $L^2(\mu)$ defined from the 
Hutchinson measure $\mu$, see (2.3). 

Both of our main results will have the {\it total\/} property for $K^\circ$ as an 
assumption. The way to test it in applications is to rely on our earlier paper \cite{JP2} 
about {\it spectral pairs\/}, i.e., subsets $\Om$, and $\Lam$, in $\br^n$ such that $\Om$ 
has finite positive $n$-dimensional Lebesgue measure, and the exponentials $\{e_\lam : 
\lam \in \Lam\}$ form an orthogonal basis for $L^2(\Om)$. We show that, for every such 
pair, the set 
$$
K:=\Lam^\circ = \{\xi \in \br^n:\xi \cdot \lam \in \bz,\qquad \forall \lam \in \Lam\}
$$
is a {\it lattice}. Analogously to the situation in Lemma 4.1 above, we also show in 
\cite{JP2} that the set $\Om$ in a spectral pair embeds in the torus $\br^n/K$. We identify 
a special class of spectral pairs, called {\it simple factors\/} which produce two finite sets 
$B$, $L\subset \br^n$, and a matrix $R$ with $K \in \lat(R)$ such that the system 
$s=(R,B,L,K)$ satisfies the conditions from section 4 above. Our present paper is 
motivated by getting ``invariants'' for simple factors from iteration of the affine maps (see 
(2.1) and (4.1) above). In \cite{JP6} we further study the converse problem of 
reconstructing simple factors from ``fractal'' iteration limit-objects. In any case, the fractal 
limit $\X(B)$ from (2.4) will also be embedded in the torus $\br^n/K$. When equipped 
with Haar-measure $L^2(\br^n/K)$ has the exponentials $\{e_\lam :\lam \in K^\circ\}$ as 
an orthogonal basis. In testing for our totality condition relative to $L^2(\mu)$, we can 
then use that $\X(B)$ is the support of $\mu$, and then apply Stone-Weierstrass to 
$\{e_\lam \}_{\lam \in K^\circ }$ when viewed as a subset of $C(\X(B))$.

We shall say that $\cL(L)$ is {\it maximal\/} if $\{e_\lam : \lam \in \cL(L)\}$ is {\it 
orthogonal\/} in $L^2(\mu)$ and (considering $t \in \br^n$) if 
$$
\text{\it whenever}\quad\langle e_t, e_\lam \rangle_{\mu} = \hm (\lam - t) =0\text{ \it for 
all } \lam \in \cL(L),\text{\it then\/ }t \in \cL(L). \tag 5.4
$$
We have used the {\it transform\/} $\hm$ given by
$$
\hm(s)=\int e_s\,d\mu=\int e^{i2\pi s \cdot x} \,d\mu(x)\qquad\text{for}\quad s \in \br^n.
\tag 5.5
$$
We say that the system $s$ is $\Lam$-orthogonal, if the functions $\{e_\lam : \lam \in 
\Lam\}$ are orthogonal in $L^2(\mu)$, here $\Lam$ is given by (4.10), i.e., 
$\Lam=L+R^* K^\circ$. (See also (5.6) below.) \enddefinition

We are now ready for the
\proclaim{Theorem 5.2} Let $s=(R,B,L,K)$ be a selfadjoint system in $\br^n$, and 
assume\roster
\item"{(i)}" $K^\circ$ is total;
\item"{(ii)}" $\cL(L)$ is maximal in $L^2(\mu)$; and 
\item"{(iii)}" $s$ is $\Lam$-orthogonal, i.e., the points in $\Lam$ from (4.10) are 
orthogonal for $\ell\neq\ell '$ in $\Lam$.\endroster
Then it follows that $s$ is of Hadamard type; i.e., $|B|=|L|$ and the $B/L$-matrix $U$ 
satisfies (3.2). \endproclaim

\demo{Proof} Condition (i) states that the orthogonal complement of $\{e_s : s \in 
K^\circ\}$ in $L^2(\mu)$ is zero. Notice that condition (iii) is equivalent to:
$$
\hm (\ell - \ell ' +R^*s)=0,\qquad\forall \ell \neq\ell '\text{ in }L,\quad\forall s \in K^\circ. 
\tag 5.6
$$
If we set
$$
\cB(t):=|B|^{-1} \sum_{b \in B}\langle b,t \rangle,\qquad \forall t \in \br^n;
\tag 5.7
$$
then (2.4) implies the factorization:
$$
\hm(t) = \cB(t) \hm (R^{*^{-1}}t). \tag 5.8
$$

For distinct points $\ell$ and $\ell '$ in $L$ we claim that
$$
R^{*^{-1}}(\ell - \ell ')\notin K^\circ.
$$
Assuming the contrary, there would be some $s \in K^\circ$ such that $\ell - \ell '=R^* s$. 
From the orthogonality property (5.6) (see Definition 5.1), we then get 
$$
\hm(0)=\hm (\ell - \ell ' -R^* s)= \langle e_\ell, e_{\ell ' + R^{*}s}\rangle _\mu =0
$$
contradicting $\hm(0)=1$.

Since $\cL (L) \subset K^\circ$, the point $t:=R^{*^{-1}}(\ell - \ell ')$ is not in $\cL(L)$. 
From the maximal property (5.4), we conclude that there is some $\lam \in \cL(L)$ such 
that 
$$
\langle e_t,e_\lam \rangle _\mu \neq 0.
$$
This term works out to
$$
\hm(R^{*^{-1}}(\ell - \ell ')-\lam)(\neq 0).
$$
From the orthogonality (Definition 5.1), we also have:
$$
0=\hm(\ell -\ell ' -R^* \lam) = \cB (\ell - \ell ')\hm (R^{*^{-1}}(\ell -\ell ')-\lam)
$$
where the last factor is non-zero. It then follows that $\cB(\ell - \ell ')=0$.

Recall, for $u,v \in \br^n$, the notation $\langle u,v \rangle :=e^{i2 \pi u \cdot v}$. Then 
the vectors $\{\langle \cdot, \ell \rangle\}$ are indexed by points $\ell \in L$, and we 
showed that they are orthogonal when viewed as elements in $\ell ^2(B)$.

It follows that $|L|\leq |B|$ where the symbol $|\cdot|$ denotes cardinality.  We claim that 
they are equal. For suppose the contrary, viz., $|L| <|B|$. Then pick coefficients $k_b \in 
\bc$, not-all zero, indexed by $b \in B$, such that
$$
\sum_{b \in B} k_b \langle b, \ell \rangle =0,\qquad \forall \ell \in L. \tag 5.9
$$
For every $s \in K^\circ$ and $\ell \in L$, consider $t:= \ell +R^* s$; and define
$$
f:=\sum_{b \in B} \bar k _b \chi_{(b+R^{-1}(\X))} \tag 5.10
$$
where $\chi$ denotes ``indicator function'', the subscript is a $b$-translate, and finally 
$\X$ is the $B$-fractal. (Recall, details below, it is compact, and satisfies $\X=B+R^{-
1}(\X)$, with
$$
\mu((b+R^{-1}(\X))\cap (b' + R^{-1}(\X)))=0\qquad\text{for all $b\neq b'$ in $B$.}) \tag 
5.11
$$
Note, (5.11) is a consequence of the totality of $K^\circ$ and the following observation: if 
$b$ and $c$ are in $B$ and $b+R^{-1}x=c+R^{-1}y$, then $R(b-c)=y-x$; it now follows 
from Lemma 4.1 that $x\in y+K$. But then
$$\align
\langle f,e_t \rangle_\mu &= \int _{\br^n} \overline{f(x)}e_t (x)\,d\mu(x)\\
&=\sum_{b \in B}k_b \int _{R^{-1}(\X)} e_t (b+x)\,d\mu(x)\\
&=\sum_{b \in B} k_b \int_{R^{-1}(\X)} \langle \ell + R^*s,b+x\rangle \, d\mu (x)\\
&=\undersetbrace (5.9) \to {\left( \sum_{b \in B} k_b \langle \ell, b\rangle \right)} 
\int_{R^{-1}(\X)}e_t \, d\mu\tag 5.12
\endalign$$
(where we use (5.11) and Lemma 5.3 below), and
$$
\langle R^* s, b \rangle = \langle s, Rb \rangle =1.
$$
The last fact is from axiom (4.6) which makes $Rb \in K$. It follows (from (5.9)) that $f$ 
is in the orthogonal complement of $\{e_{\ell +R^* s}\}$ as $\ell$ varies over $L$, and 
$s$ over $K^\circ$. But from (i), we know that this is a total set of vectors in 
$L^2(\mu)$, so the function $f$ must vanish identically, $\mu$-a.e. If the coefficients 
$\{k_b\}$ are not all zero, this would contradict (5.11), (2.3), and the basic properties of 
the Hutchinson measure $\mu$.

From the contradiction, we conclude that $|L|=|B|$; which is to say, both conditions on 
the matrix $\left( \langle b,\ell \rangle \right)_{b,\ell}$, indexed by $B\times L$, to be of 
generalized Hadamard type, are satisfied. We have $|L|=|B|=N$. If the matrix is denoted 
$U$, then 
$$
UU^*=U^*U=NI_N \tag 5.13
$$
where $I_N$ denotes the identity matrix in $N$ variables, and $U^*$ is the transpose 
conjugate. To define it, it is convenient to use a common index labeling, e.g., 
$\{1,2,\ldots,N\}$.\qed
\enddemo

\remark{Remark} Note that if we assume $B$ is a subset of a set of representations for 
$R^{-1}K/K$, then (5.11) follows from an application of \cite{Ke, Theorem 10} and a 
related result in \cite{Ma}. (See also \cite{Ba-Gr} for related work.) The Kenyon-Madych 
result applies in the present context since the mapping {\it from\/} the set of all finite $B$-
strings $(b_1,\ldots,b_m)$ with $m$ varying in $\bn$, $b_i \in B$, {\it into\/} $\sum_i R^i 
b_i \in K$ is 1-1. This follows by induction and use of our orthogonality assumptions. To 
use \cite{Ke}--\cite{Ma}, we then extend $B$ so as to get a full set of residue classes 
$R^{-1}(K)/K$. \endremark

\remark{Question} Does either of the following two conditions imply the other: (i) 
$K^\circ$ is total, (ii) $B$ is a subset of a set of representatives for the quotient $(R^{-
1}K)/K$\,?\endremark

In the calculation (5.12) above, the following lemma was used. (It is needed because we 
do {\it not\/} know if, in general, $\mu$ is a Hausdorff-measure.)

\proclaim {Lemma 5.3} Under the assumptions of Theorem 5.2, it follows that
$$
\int_{\sig_b \X} f(x)\,d\mu(x)=\int_{R^{-1}\X} f(x+b)\, d\mu(x)
$$
for all $b$ in $B$ and $f$ in $L^2(\mu)$.\endproclaim

\demo{Proof} The claim is equivalent to having
$$
\mu(\sig_b\Delta)=\mu(R^{-1}\Delta)
$$
for all $\mu$-measurable sets $\Delta\subset\br^n$ and all $b$ in $B$; which in turn is 
equivalent to
$$
\mu(\sig_b\Delta)=\mu(\sig_c\Delta) \tag 5.14
$$
for all $\mu$-measurable $\Delta$, and all $b$ and $c$ in $B$. The last equivalence used 
the assumption that $0\in B$. By regularity, it suffices to consider the case where $\Delta$ 
is a closed set.

Let $\Delta$ be a closed subset of $\X(B)=\X$, and choose $B_k\subset B^k=B\times 
\cdots \times B$ ($k$ terms) such that
$$
(b_1,\ldots,b_k,b_{k+1})\subset B_{k+1}\Rightarrow (b_2,\ldots,b_{k+1})\in B_k,
$$
and such that
$$\Delta=\bigcap_{k=1}^\infty E_k,
$$
where 
$$
E_k=\bigcup_{(b_1,\ldots,b_k)\in B_k} \sig_{b_1}\cdots\sig_{b_k}\X.
$$
The first condition means that $(E_k)$ is a descreasing sequence of compact sets. It 
follows (analogously to (5.11)) that the overlaps in the definition of $E_k$ are $\mu$-
null-sets. Hence, to prove the lemma, it suffices to show that
$$
\mu(\sig_{b_1}\cdots \sig_{b_k}\X)=\mu(\sig_{c_1}\cdots\sig_{c_k}\X)
$$
for all $(b_1,\ldots,b_k)$, $(c_1,\ldots,c_k)$ in $B^k$.

First note that, for any Borel set $\Delta$, and any $b$ in $B$, we have
$$
\mu(\Delta)=|B|^{-1}\sum_{c\in B}\mu(\sig_c^{-1}\Delta)
\geq |B|^{-1}\mu(\sig_b^{-1}\Delta),
$$
and hence $\mu(\sig_b\Delta)\geq|B|^{-1}\mu(\Delta).$ From this inequality, and (5.11), it 
follows further that
$$
\mu(\X)=\mu(\bigcup_{b \in B}\sig_b\X)=\sum_{b\in B}\mu(\sig_b\X)
\geq \mu(\X),
$$
and therefore that $\mu(\sig_b\X)=|B|^{-1}\mu(\X)$. Assuming
$$
\mu(\sig_{b_2}\cdots\sig_{b_k}\X)=|B|^{-k+1}\mu(\X),
$$
it follows (analogously to the above), that
$$
\mu(\sig_{b_1}\sig_{b_2}\cdots\sig_{b_k}\X)=|B|^{-k}\mu(\X).\tag 5.15
$$
Hence, by induction, (5.15) is true for all positive integers $k$, and all $b_1,\ldots, b_k$ 
in $B$. This completes the proof of the lemma.\qed
\enddemo

Note also (heuristically) that (5.14) is a consequence of (5.11), and that
$$\align
\mu(\sig_b\Delta)&= \lim_{k\rightarrow \infty}\sum_{b_1,\ldots,b_k}|B|^{-1}
\chi_{\sig_b \Delta}(\sig_{b_1}\cdots \sig_{b_k}x)\\
&=\lim_{k\rightarrow \infty}\sum_{b_1,\ldots,b_k}|B|^{-1}
\chi_\Delta((\sig_{b_2}\cdots\sig_{b_k}x)+Rb_1 -Rb)
\endalign$$
where $x$ in $\X$ is arbitrary. However, the first equality requires that $\chi_{\sig_b 
\Delta}$ is continuous. We refer to \cite{Fa, p. 121} for further details on this point.

We conclude with the following lemma which is both basic and general; in fact it holds in a 
context which is more general than where we need it. Such more general contexts occur, 
e.g., in \cite{St4}, \cite{Mat}, \cite{MOW}, and \cite{Od}, (among other places). But we 
will still restrict the setting presently to where it is needed below for our proof of Theorem 
6.1.

\proclaim{Lemma 5.4} Let $(R,B)$ be an affine system in $\br^n$ (see details in Section 
2) with $R$ expansive and $B \subset \br^n$ a finite subset. Let $\cB$ be given by 
\rom{(5.7)}, and let $\mu$ be the probability measure from \rom{(2.2)}. We are assuming 
the property \rom{(5.11)}. Let $\Cal N:=\{t \in \br^n :\cB(t)=0\}$. Then, for the roots of 
$\hm$, we have
$$
\{t \in \br^n:\hm(t)=0\}=\bigcup_{k=0}^\infty R^{*^k}(\Cal N).
$$
\endproclaim

\demo{Proof} We have (5.8) by virtue of \cite{JP6, Lemma 3.4}, and it follows that 
$\hm(t)=0$  when $t \in R^{*^k}(\Cal N)$ for some $k \in \{0,1,\ldots\}$. From the 
assumed expansivity of $R$, we also know that the corresponding infinite product formula 
is convergent, see \cite{JP6, (3.13)}. In fact $\lim_{k \rightarrow \infty}\hm (R^{*^{-
k}}t)=1$, for all $t \in \br^n$. This is from continuity of $\hm$, and the limit, $R^{*^{-
k}}t \rightarrow 0$. Now consider,
$$
\hm(t)=\prod_{j=0}^{k-1} \cB(R^{*^{-j}}t) \hm(R^{*^{-k}}t),
$$
and suppose $\hm(t)=0$. Pick $k$ (sufficiently large) s.t., $\hm(R^{*^{-k}}t)\neq 0$. 
(This is possible by continuity, and the fact that $\hm(0)=1$). We conclude, then that, for 
some $j$, $0\leq j <k$, $R^{*^{-j}}t\in \Cal N$; and this is the assertion of the 
lemma.\qed\enddemo

\head 6. Orthogonal Exponentials \endhead
We keep the standing assumptions on the quadruple $s=(R,B,L,K)$ which determine a 
system in $\br^n$. In particular, the matrix $R$ is assumed {\it expansive\/} (see section 
2), the sets $B$ and $L$ in $\br^n$ are finite both containing $0$. We will assume now 
that $s$ is selfadjoint and of Hadamard type. We say that the system $s$ is {\it 
irreducible\/} if there is no proper linear subspace $V\subset \br^n$ (i.e., of smaller 
dimension) which contains the set $B$, and which is invariant under $R$, i.e., $Rv \in V$ 
for all $v\in V$. If such a proper subspace does exist, we say that $s$ is {\it reducible}. In 
that case, it is immediate that the fractal $\X(B)$ from (2.4) is then contained in $V$. All 
the examples in Section 7 below can easily be checked to be irreducible. But the following 
example in $\br^2$ is reducible, and serves to illustrate the last conclusion from our 
theorem in the present section: Let $R=\pmatrix 2&1\\0&2\endpmatrix$, 
$B=\left\{\pmatrix0\\0\endpmatrix, \pmatrix\frac{1}{2}\\0\endpmatrix\right\}$, $L=\left\{ 
\pmatrix 0\\0\endpmatrix, \pmatrix 1\\0 \endpmatrix \right\}$, $K=\bz^2$; and 
$V=\pmatrix \br\\0\endpmatrix$, i.e., the $x$-axis in $\br^2$. Then a direct calculation 
shows that $\X(B)=\pmatrix I\\0 \endpmatrix$, where $I=[0,1]$ is the unit-interval on the 
$x$-axis; and the Hutchinson measure $\mu$ is the product measure $\lam_1 \otimes 
\delta_0$ where $\lam_1$ is the restriction to $I$ of the one-dimensional Lebesgue 
measure, and $\delta_0$ is the point-$0$-Dirac measure in the second coordinate. For the 
set $\cL(L)$ from (5.1), we have,
$$
\cL(L)=\left\{ \pmatrix n\\N(n)\endpmatrix : n=0,1,2,\ldots \right\},
$$
$n$ represented by finite sums, $n=\sum_{j \geq 0}2^j\epsilon_j$, $\epsilon_j \in \{0,1\}$, 
and $N(n)=\sum_{j>0} j2^{j-1}\epsilon_j$.

We will show below that if $K^\circ$ is total then in the ``{\it fractal case\/}'', i.e., when 
$N=|B|=|L|$ is smaller than $|\det R|$, and when $s$ is {\it irreducible}, then $\cL(L)$ is 
{\it maximally orthogonal}, and $s$ is $\Lam$-orthogonal. Recall $\mu$ is the 
Hutchinson measure, see (2.3), $\cL(L)$ is the set given by (5.1); and finally the 
properties regarding the two sets $K^\circ$ and $\cL(L)$ refer to the corresponding 
exponentials $e_s$, when $s$ is in the respective sets, and each
$$
e_s(x)=\langle s,x \rangle =e^{i2 \pi s \cdot x}\tag 6.1
$$
is considered a vector (alias function on $\br^n$) in the Hilbert space 
$L^2(\mu)=L^2(\X,\mu)$. We refer to Definition 5.1 and Theorem 5.2 for further details. 
Recall that the total property (iii) in Theorem 5.2 amounts to the $\Lam$-orthogonality, 
including the assertion
$$
\hm(\ell - \ell ' + R^*s)=0 \tag 6.2
$$
for all $\ell \neq \ell '$ in $L$ and all $s \in K^\circ$. But if $s$ is selfadjoint and of 
Hadamard type, then (6.2) follows immediately from (5.8), which is the functional 
equation of the transform $\hm$, see also (5.5).

The purpose of the present section is twofold. First we show that Theorem 5.2 has a 
partial converse, and secondly that the technical conditions from our two theorems 5.2 and 
6.1 amount to the dual pair condition (see Section 1) for representations of the $C^*$-
algebra $\cO_N$. This is for systems $s=(R,B,L,K)$ as specified where the two given 
finite sets $B$ and $L$ in $\br^n$ are assumed to have the same cardinality $N$, i.e., 
$|B|=|L|=N$. Our recent paper \cite{JP6} further details how the representation duality 
relates to our present assumptions. But we shall summarize the essentials here for the 
convenience of the reader. The Cuntz-algebra $\cO_N$ (see \cite{Cu}) is known to be 
given universally on $N$ generators $\{s_i\}$ and subject only to the relations: 
$$s_i^*s_j = \delta_{ij}1\andd \sum_i s_is_i^*=1 \tag 6.3 $$

\noindent This means that, if a finite set of $N$ operators $S_i$ say, acting on some 
Hilbert space $\cH$ say, are known to satisfy the relations (6.3), then there is a unique 
represention $\rho$ of $\cO_N$, acting by bounded operators on $\cH$, such that 
$\rho(s_i)=S_i$ for all $i$; or, equivalently, $\rho(a)f=\hat a f$ for all $a \in \cO_N$ and 
all $f\in \cH$, where the operator $\hat a$ is given by the same expression in the $S_i$s as 
$a$ is in the $s_i$-generators.

For a given system $s=(R,B,L,K)$, there is then the possibility of making a representation 
duality based on the exponentials $e^{i2\pi t\cdot x}$ in (6.1),  and treating the two 
vector-variables $x$ and $t$ symmetrically: The pair $(R,B)$ gives one affine system 
$\sigma_b x=R^{-1}x+b$ ($b\in B$) in the $x$-variable; and the dual system $(R^* ,L)$ 
given by $\tau_\ell t=R^*t+\ell$ ($\ell \in L$); a second one now acting in the $t$-variable. 
See (2.1) and (4.1) above. To be able to generate the asserted representation pair we need 
to specify $\{\tau_\ell\}_{\ell \in L}$ for enough values of $t$ such that the corresponding 
functions $e_t(x):=e^{i2\pi t\cdot x}$ span a dense subspace in $L^2(\mu)$. But it turns 
out that other conditions must be met as well: For the two affine systems 
$\{\sig_b\}_{b\in B}$ and $\{\tau_\ell \}_{\ell\in L}$, the question is if we can associate 
operator systems $\{S_b\}$ and $\{T_\ell\}$ of $2N$ operators acting on $L^2(\mu)$, 
each system satisfying (6.3), and the operators collectively defined from the exponentials 
$e^{i2\pi t \cdot x}$ as specified. When this is so, we have an (orthogonal) dual pair of 
representations of $\cO_N$ acting on $L^2(\mu)$, and conversely. Then it turns out that, 
for each pair $(b,\ell)$, the operator $S^*_bT_\ell$ is a multiplication operator, see (6.6) 
below.

\proclaim {Theorem 6.1} Let $s=(R,B,L,K)$ be a system in $\br^n$ and assume that $s$ 
is selfadjoint and of Hadamard type. Assume further that $K^\circ$ is total (with a minimal 
choice for $K$), and that $|B|<|\det R|$, and let $\mu$ be the corresponding measure (see 
(2.3)) with support $\X$. Then $s$ is $\Lam$-orthogonal and carries a dual pair of Cuntz 
representations (with $\cO_N$ acting on $L^2(\mu)$ for both representations). If $s$ is 
also irreducible, then $\cL(L)$ is maximally orthogonal. \endproclaim

\demo{Proof} Let $N=|B|=|L|$ and note that from \cite{JP6} (Theorem 4.1) we get a 
dual pair of representations $\{S_b\}_{b\in B}$ and $\{T_\ell \}_{\ell \in L}$ of the 
Cuntz algebra $\cO_N$, see also \cite{Cu} and \cite{Ar}, acting on $L^2(\mu)$ and 
given by the respective formulas:
$$\align
S_b^*f&= N^{-1/2}f\circ \sig_b\qquad\text{for}\quad f\in L^2(\mu),\quad b\in B
\tag 6.4\\
T_\ell e_s&=e_{\tau_{\ell}(s)}\qquad \text{for}\quad s \in K^\circ\quad \text{and}\quad 
\ell \in L. \tag 6.5
\endalign$$

Moreover $S^*_b T_\ell$ is the multiplication operator $M_{b\ell}$ on $L^2(\mu)$ given 
by
$$
M_{b\ell}f=N^{-1/2}(e_\ell \circ \sig_b)f\qquad\text{for}\quad f \in L^2(\mu).
\tag 6.6
$$
It follows from \cite{JP6, Theorem 4.1} that for $\ell$ and $\ell '$ in $L$
$$
T^*_\ell T_{\ell '}=\delta_{\ell\ell '}I\andd \sum_{\ell \in L}T_\ell T_\ell^*=I 
\tag 6.7
$$
Here we use the Kronecker delta notation
$$
\delta_{\ell \ell '}=\cases 1&\text{if $\ell=\ell '$}\\
      0&\text{if $\ell\neq \ell '$},\endcases
$$
and $I$ denotes the identity operator in the Hilbert space $L^2(\mu)$.

It follows then from (6.6) that the vectors $e_{\ell +R^*s}$, $s \in K^\circ$, are mutually 
orthogonal in $L^2(\mu)$ for distinct values of $\ell$, i.e., for $\ell\neq\ell '$ in $L$. For 
more details on this point, we refer to sections 3--4 in \cite{JP6}. For 
$$
\lam =\sum_{j=0}^n R^{*^j}\ell_j, \tag 6.8
$$
we have
$$
T_{\ell_0}T_{\ell_1}\cdots T_{\ell_n}e_0=e_\lam. \tag 6.9
$$

Let $\lam=\sum_{j=0}^n R^{*^j}\ell_j$, and $\kappa=\sum_{j=0}^mR^{*^j}k_j$, where 
the $\ell_j$'s and $k_j$'s are in $L$. Then $e_\lam$ and $e_\kappa$ are orthongonal in 
$L^2(\mu)$ except in the cases where $m \leq n$, and $\ell_j=0$ for $j>m$, and where 
$m\geq n$ and $k_j=0$ for $j>n$. In the exceptional cases, it follows from (6.9) and $T_0 
e_0=e_0$ that $e_\lam=e_\kappa$. To prove the orthogonality assertion above, note that
$$
\langle e_{\lam},e_\kappa \rangle_\mu =
\langle T_{\ell_0}T_{\ell_1}\cdots T_{\ell_n} e_0,
T_{k_0}T_{k_1}\cdots T_{k_m}e_0\rangle
$$
is $=0$, unless $\ell_0=k_0$, because $T_{k_0}^*T_{\ell_0}=0$ if $k_0\neq \ell_0$. If 
$k_0=\ell_0$, then $T_{k_0}^*T_{\ell_0}=I$, and we can repeat the argument on 
$\ell_1$ and $k_1$. It remains to consider the case where $n>0$ and $m=0$; in this case, 
we will use the identity, $T_0 e_0=e_0$, to write
$$
\langle e_\lambda,e_\kappa\rangle_\mu=
\langle T_{\ell_0} T_{\ell_1} \cdots T_{\ell _n} e_0,T_0 e_0\rangle_\mu=0
$$
and we conclude that $\langle e_\lam,e_\kappa \rangle_\mu =0$ unless 
$\ell_0=\ell_1=\cdots=\ell_n=0$.

It follows that the map,
$$
(\ell_0,\ldots,\ell_n)\mapsto \sum_{j=0}^n R^{*^j}\ell_j \in \cL(L)\tag 6.10
$$
is $1$--$1$  on the set of finite sequences $(\ell_0,\ldots , \ell_{n-1},\ell_n)$ with $n$ a 
nonnegative integer, the $\ell_j$'s in $L$, and $\ell_n\neq 0$. 

The assumption that $s$ be irreducible is now imposed, and we show that $\cL(L)$ has 
the stated maximality property: We show that, if $t \in \br^n$ and $\langle e_\lam, e_t 
\rangle_\mu =0$ for all $\lam \in \cL(L)$, then it follows that $t \in \cL(L)$. We shall do 
this by contradiction, assuming the $t \notin \cL(L)$. We shall use the functional equation 
(5.8) for $\hm$, recalling that
$$
\langle e_\lam,e_t \rangle_\mu=\hm(t-\lam). \tag 6.11
$$
We shall also use that for every $s \in \br^n$ there is some $\ell \in L$ such that $\cB(\ell -
s)\neq 0$. This follows from the formula (5.7) for $\cB(\cdot)$, and from the Hadamard 
property (3.2) which is now assumed.

As a special case of (5.8), we get
$$
0=\hm (t-\ell_0-R^*\ell_1)=\cB (t-\ell_0)\hm (R^{*^{-1}}t-R^{*^{-1}}\ell_0 -\ell_1).
$$
Picking $\ell_0 \in L$ s.t. $\cB(t-\ell_0)\neq 0$, we get
$$\align
0&=\hm(R^{*^{-1}}t - R^{*^{-1}}\ell_0 - \ell_1 - R^*\ell_2)\\
&=\cB (R^{*^{-1}}t-R^{*^{-1}}\ell_0 - \ell_1)
\hm(R^{*^{-2}}t- R^{*^{-2}}\ell_0 - R^{*^{-1}}\ell_1 - \ell_2).
\endalign$$
Picking $\ell_1 \in L$ s.t.
$$
\cB(R^{*^{-1}}t - R^{*^{-1}}\ell_0 - \ell_1)\neq 0,
$$
we conclude next that
$$
\hm(R^{*^{-2}}t - R^{*^{-2}}\ell_0 - R^{*^{-1}}\ell_1 - \ell_2) =0,
$$
and we continue by induction, determining $\ell_0,\ell_1,\ldots \in L$ such that the points 
$$
s_p := R^{*^{-p}}\ell_0+\cdots +R^{*^{-1}}\ell_{p-1}+\ell_p
$$
are in the dual fractal set $\X(L)$, see (6.2) above. When $N<|\det R|$, we may pick, 
inductively, the ``digits'' $\ell_i$ such that the differences
$$
R^{*^{-p}}t-s_p \tag 6.12
$$
are distinct as $p$ varies, but
$$
\hm(R^{*^{-p}}t-s_p)=0 \andd \cB(R^{*^{-(p-1)}}t - s_{p-1}) \neq 0.
$$
Notice that the analytically extended transform 
$$
\hm(z)=\int e^{i2 \pi z \cdot x}\,d\mu(x)\tag 6.13
$$
is {\it entire\/} analytic on $\bc^n$, where for $z=(z_1,\ldots,z_n) \in \bc^n$, $z \cdot x = 
z_1x_1 + \cdots + z_n x_n$ is the usual dot-product. Hence its zeros cannot accumulate. 
But the ``dual attractor'' $\X(L)$ (see (7.2)) is compact in $\br^n$ so there a subsequence 
$s_{p_i} $ with limit $s_{p_i}\rightarrow s \in \X(L)$, and 
$$
0=\lim_{p_i} \hm(R^{*^{-p_i}}(t) - s_{p_i})=\hm(-s) = \overline{\hm(s)}
$$
contradicting that the roots of $\hm(\cdot)$ must be isolated (see (6.12)), even isolated in 
$\bc^n$. The contradiction completes the proof, and we conclude that $\cL(L)$ is 
maximal.\enddemo

If only a finite number of the ``digits'' $\ell_j$ are nonzero, then, using the contractive 
property of $R^{*^{-1}}$, we see that the sequences 
$R^{*^{-p}}(t)$, and $s_p$, both converge to zero as $p \rightarrow \infty$, 
contradicting that $\hm(0)=1$, since $\lam \rightarrow \hm (\lam)$ is continuous on 
$\br^n$.

\definition{Claim 1} The set $B$ is a subset of a set of representatives for $R^{-1} 
(K)/K$. \enddefinition

\demo{Proof of Claim} From the self-adjointness of $s$ we have $RB \subset K$ (by 
Lemma 4.2). Therefore, $B \subset R^{-1}K$. If $b$ and $b'$ are distinct and both in 
$B$, and if $b \in b' +K$, then $e_t (b+x)=e_t (b' +x)$ for all $x \in R^{-1}\X$ (all $t \in 
K^\circ$), contradicting the totality of $K^\circ$ in $L^2(\mu)$.\enddemo

\definition {Claim 2} The finite set $L$ is a subset of a set of representatives for 
$K^\circ/R^*K^\circ$.\enddefinition

\demo{Proof of Claim} By Lemma 4.2, $L\subset K^\circ$. If $\ell$ and $\ell '$ are in 
$L$, and $\ell=\ell '+R^* \gamma$ for some $\gamma\in K^\circ$, then 
$$\align
\langle b,\ell \rangle &=\langle b, \ell ' +R^*\gamma\rangle\\
&=\langle b, \ell ' \rangle \langle Rb,\gamma \rangle\\
&=\langle b,\ell ' \rangle
\endalign$$
where the last equality used Lemma 4.2 again. But this contradicts the Hadamard-
property, unless $\gamma=0$. Considering,
$$\align
x_p&=R^{*^{-p}}t-(R^{*^{-p}}\ell_0+\cdots +R^{*^{-1}}\ell_{p-1} +\ell_p)
\, ,\quad\text{and}\\
y_p&=R^{*^p}x_p=t-(\ell_0+R^*\ell_1 + \cdots + R^{*^p}\ell_p )\, ,
\endalign$$
and letting $P=\{p:\ell_p \neq 0\}$; then we showed above that $P$ is infinite, and that $p 
\in P\rightarrow y_p$ is a $1$--$1$ map. Hence $\{y_p :p\in P\}$ is infinite.\qed 
\enddemo

\remark{Remark \rom{6.2}} For the {\it reducible\/} example (in $\br^2$) mentioned in 
the beginning of the present section, we note that all the conditions of the first part of 
Theorem 6.1 are satisfied. We also described the set $\cL(L)$ of orthogonal exponentials 
for the example. But the maximality condition is not satisfied relative to $L^2(\mu)$. 
Indeed, for the transform $\hm(s)$ from (5.5), we have, with $s=(s_1, s_2)\in \br^2$,
$$
\hm(s)=\left\{ \matrix e^{is_1 \pi}\, \frac{\sin(s_1 \pi)}{s_1 \pi} &\text{if}&(s_1\neq 0)\\
1 & \text{if}&s_1=0.\endmatrix\right.\
$$
It follows that the identity from (5.4) will be satisfied whenever  $t=(t_1,t_2)\in \br^2$ is 
such that $t_1\in\bz_{-}$, i.e., negative and integral. (Specifically, $\hm(\lam - t)=0$ for 
$\forall \lam \in \cL(L)$.) From the calculation of $\cL(L)$, we note that such points 
$t=(t_1,t_2)$ will {\it not\/} be in the set $\cL(L)$; and so the maximality condition is not 
satisfied for the example.\endremark

\head 7. Examples \endhead
\subhead 7.1 Background Material \endsubhead
We now give examples to illustrate the conditions in Theorems 5.2 and 6.1. Since the 
generalized Hadamard matrices are known up to $N=4$, the examples we give are 
``typical'' for the possibilities when $N\leq 4$, and it is likely that there is a {\it 
classification\/}; but as it is unclear what is the ``correct'' notion of {\it equivalence\/} for 
systems $s=(R,B,L,K)$ we will postpone the classification issue to a later paper. Note that 
the examples occur in pairs, one for $s$ and a {\it dual one\/} for $s^\circ$. Also note that 
each $s$ will correspond to a {\it spectral pair\/} $(\Om,\Lam)$ as well as a selfsimilar {\it 
iteration limit\/}, typically a ``fractal'' $\X$ with a {\it selfsimilar\/} measure $\mu$. When 
the given system $s$ is {\it selfadjoint\/}, then there will in fact be a {\it pair\/} of 
``fractals'' occurring as iteration limits, a selfsimilar $\mu$ from the affine system:
$$
\sig_b x=R^{-1}x+b\qquad\text{leading to $\X=\X(B)$}\tag 7.1
$$
defined from $s$, and also
$$
\tau^{-1}_{\ell}(t)=R^{*^{-1}} (t-\ell)\qquad\text{leading to $\X(L)$}, \tag 7.2
$$
and defining the corresponding dual selfsimilar measure $\mu '$. Recall both $\mu$ and 
$\mu '$ are probability measures on $\br^n$; $\mu$ is determined by (2.2), and $\mu '$ 
by:
$$
\mu '=|L|^{-1} \sum_{\ell \in L}\mu ' \circ \tau_\ell, \tag 7.3
$$
see also (7.1)--(7.2) and Lemma 4.1 for more details on the dual pair of affine systems.

Our examples below will be constructed from the matrices (3.3)--(3.5) which we listed in 
section 3. In fact, we shall supply a group of examples for each of the generalized 
Hadamard matrices $N=2$, $N=3$, and $N=4$, all the examples will be symmetric and of 
Hadamard type; but some will {\it not\/} be selfadjoint. In fact, when considering 
$s=(R,B,L,K)$ we shall fix the first three $R$, $B$, and $L$, but allow variations in the 
lattice. When we insist on the Hadamard type, we shall see that, in some familiar fractal-
examples, it will then {\it not\/} be possible to choose {\it any\/} lattice $K$ such that the 
corresponding system $s=s(-,K)$ is selfadjoint. We will then say that the system is {\it 
not\/} self-adjoint; it turns out that the obstruction is a certain {\it integrality condition\/}; 
and, when it is not possible to find a lattice consistent with both selfadjointness and 
Hadamard type, then it will typically be a simple, case by case computation, and we shall 
be very brief with detailed calculations. (It will be immediate that each of the examples in 
the list is {\it irreducible}; see section~6.)

\subhead 7.2 Group 1 Examples \endsubhead
We take $N=2$; the matrix is (3.3), and the examples are illustrated with subsets of the 
line, i.e., $n=1$, for $\br^n$. First, take $R=4$, i.e., multiplication by the integer 4; the 
sets $B$ and $L$ will be $B=\{0, 1/2\}$, $L=\{0,1\}$, and lattice $K=\bz$. The 
$(\Om,\Lam)$ spectral pair will be as follows:\roster
\item"{$\Om =$}" $[0,1/4]\cup [1/2,3/4]$, (i.e., the union of two intervals);\medskip 
\item"{$\Lam =$}" $\{0,1\}+4\bz $, (i.e., two residue systems modulo 4, see (4.10) above 
for the general case) \medskip 
\item"{$\X =$}" iteration fractal, see Figure A, fractal dimension $D={\ln 2 \over \ln 
4}=1/2$, see the affine system (2.2), and also more details on $\mu$ in Section 2 of 
\cite{JP6}.\endroster\par

It is easy to check that with this choice for $R$, $B$, $L$, and $K$, the corresponding 
system $s$ is selfadjoint and of the Hadamard type. For this particular example, there are 
only two choices for $K$ such that the corresponding system $s_K =s(-,K)$ is selfadjoint. 
They are $K=\bz$ and $K=2\bz$. But the following modification, corresponding to the 
classical middle-third-Cantor set, will only be symmetric; not selfadjoint: With
$$\align
R&=3,\\
B&=\{0,2/3\},\andd\\
L&=\{0,3/4\},
\endalign$$
we have the Hadamard type, c.f., (3.2); but there is {\it no\/} lattice $K$ in $\br$ which 
makes the corresponding system $s_K$ selfadjoint (Graphic illustration, Figures A and 1).

\subhead 7.3 Group 2 Examples \endsubhead
We take $N=3$; the matrix is (3.4), and the examples are illustrated with subsets of the 
plane $\br^2$. Take
$$\align
R&={\pmatrix 6&0\\0&6 \endpmatrix},\\
B&=\left\{\binom 00,\binom{\frac12}{0},\binom {0}{\frac12} \right\},\\
L&=\left\{\binom 00,\pm \ell\right\} \qquad \text{where}\quad\ell={2 \over 
3}\binom{1}{-1},\\
K&=3\bz^2,\text{ i.e., multiples of the unit-lattice in 2 dimensions, equivalently }\\
&\qquad\text{points in $\br^2$ of the form $\binom{3m}{3n}$ where }m,n\in \bz.
\endalign$$

The corresponding system will be selfadjoint of Hadamard type. If $K$ is instead taken to 
be the lattice generated by the two vectors $\binom 11$ and $\binom 0{\frac 32}$ (which 
turns out to yield $K\subset L^\circ$), then there is a corresponding spectral pair 
$(\Om,\Lam)$ where $\Om$ is a suitable union of scaled squares in the plane, and the 
corresponding spectrum satisfies $\Lam^\circ=K$. But with this $K$, the iteration system 
$s_K$ will not have $K^\circ$ total in $L^2(\mu)$. In all, there are only three distinct 
choices, in this case, for lattices $K$ in $\br^2$ such that the corresponding system $s_K 
=s(-,K)$ is selfadjoint: They are given by the respective matrices $3I_2$, $\pmatrix 
3&0\\3&{\frac 32} \endpmatrix$, and $\pmatrix 1&0\\1&{\frac 32}\endpmatrix$ with 
inclusions $K_1 \subset K_2 \subset K_3$ for the lattices. The fractal dimension of $\X$ 
is 
$$
D={\frac {\ln 3}{\ln 6}}\simeq .61.
$$
(Graphic illustrations, Figures 2--9.)

\subhead 7.4 Group 3 Examples \endsubhead
We take $N=4$; the matrix is (3.5) corresponding to $u=-1$, and the examples are 
illustrated with solid sets, i.e., pictures in 3-space $\br^3$. Take
$$\align
R&={\pmatrix 2&0&0\\0&2&0\\0&0&2\endpmatrix},\\
B&=\left\{ {\pmatrix 0\\0\\0 \endpmatrix},{\pmatrix {-\frac 12}\\0\\0\endpmatrix},
        {\pmatrix 0\\{-\frac 12}\\0\endpmatrix},
{\pmatrix 0\\0\\{-\frac12}\endpmatrix}
       \right\},\\
L&=\left\{ {\pmatrix 0\\0\\0 \endpmatrix},{\pmatrix -1\\-1\\0 \endpmatrix},
       {\pmatrix -1\\0\\-1 \endpmatrix},{\pmatrix 0\\-1\\-1 \endpmatrix} \right\},\\
K&=\bz^3=K^\circ \qquad\text{(i.e., selfduality)}.
\endalign$$

It is convenient to summarize the choices for $R$, $B$, $L$ and $K$ as follows: $$\align
R&=2I_3,\\
B&={-\frac 12} I_3,\qquad\text{and}\\
L&={-\pmatrix 1&1&0\\1&0&1\\0&1&1 \endpmatrix},
\endalign$$
where
$$
I_3 = \pmatrix 1&0&0\\0&1&0\\0&0&1 \endpmatrix
$$
is the unit-matrix. The choices for the lattice $K$ are subjected to the conditions in 
Lemma 4.1. It turns out that the choice $K=\bz^3$ is {\it the\/} minimal one such that the 
system $s_K=s(-,K)$ is {\it selfadjoint\/}; and there is also a unique maximal choice for 
$K$ with $s_K$ selfadjoint, viz., $K=L^\circ$ where $L^\circ$ is given by (7.4) below. 
(Since $L$ is symmetric, the matrix for the lattice $L^\circ$ is $L^{-1}$.) The lattice 
$L^\circ$ has matrix represented by the inverse 
$$
L^{-1}={-\frac 12 }\pmatrix 1&1&-1\\1&-1&1\\-1&1&1\endpmatrix. \tag 7.4
$$

The corresponding system $s=(R,B,L,K)$ is selfadjoint of Hadamard type. If the choice 
for $K=\bz^3$ is replaced by $K=L^\circ$, then (the modified) $s$ is still selfadjoint: 
Notice that $K=\bz^3$ is the minimal choice for $K$ (subject to (i)--(iii) in Lemma 4.1); 
and $K=L^\circ$ is the maximal one. This means that $K^\circ = \bz^3$ is maximal among 
the possible choices for $K^\circ$; and this $K^\circ$ is total, see Definition 5.1.

\subhead 7.5 Dual Pairs \endsubhead
The fractal dimension is $D={\frac{\ln 4}{\ln2}}=2$ which is integral, but less than the 
dimension (viz., 3) of ambient $\br^3$. The fractal for the system $s$ arises from scaling 
iteration of the set $\Om =$ union of 4 cubes, see the figure (Figure B). For the dual 
system, $\Om^\circ$ is instead the union of tetrahedra resulting in a 3-dimensional 
Sierpinski gasket, same fractal dimension $D=2$, but with angles $60^\circ$ rather than 
$90^\circ$. The sketch is Figures 10--17, see also \cite{Sch} for more details; it is the 
Eiffel tower construction, (maximal strength with least use of iron.)

The corresponding planar Sierpinski-gasket corresponding to $R=\pmatrix 2&0\\0&2 
\endpmatrix$, $B={\frac 12}I_2$, and $L={\frac 23}
\pmatrix1&-1\\-1&1\endpmatrix$, does {\it not\/} have a lattice choice for $K$ which 
makes the associated system $s=(R,B,L,K)$ in $\br^2$ selfadjoint. The fractal dimension 
is $D={\frac{\ln 3}{\ln2}}\approx 1.58.$

For the matrix (3.5) with primitive 4th roots of 1 (e.g., $u=i$), there is also a realization in 
$\br^3$: We may take $R=2I_3$, $B={\frac 12}I_3$, and $L=\pmatrix {\frac 
12}&1&{\frac32}\\1&0&1\\{\frac 32}&1&{\frac 12}\endpmatrix$ will give a system of 
Hadamard type in 3-space, but again there is no lattice choice for $K$ in $\br^3$ such that 
the corresponding $s_K$ is selfadjoint. (Graphic illustrations, Figures 10--17.)

\head 8. Concluding Remarks \endhead

The operators $\{T_\ell \}$ from (6.5) and (6.7) may also be used in the definition of an 
endomorphism $\theta$ on a certain $C^*$-algebraic $\cO_N$-crossed product, $\frak 
U$ say. It is given by,
$$
\theta(A)=\sum_{\ell \in L} T_\ell A T_\ell ^*,\quad \text{for }A\in \frak U\, ,
$$
and clearly, $\theta(A^*)=\theta(A)^*$, and $\theta(AB)=\theta(A)\theta(B)$ for all $A,B 
\in \frak U$. Continuous versions, also called endomorphism-semigroups, have been 
studied recently by Arveson and Powers, see e.g., \cite{Ar}. As spectral-invariants for 
these, Arveson has proposed (in \cite{Ar})  a Cuntz-algebra construction which is based 
on Wiener-Hopf techniques, and which is inherently {\it continuous\/}, in fact with 
$\br_+$ used as index for the generators in place of the usual finite (or infinite) {\it 
discrete\/} labeling set $\{1,\ldots,N\}$. For our present $B/L$ duality project with dual 
fractals, $\X(B)$ and $\X(L)$; we plan (in a sequel paper) to study an analogous $C^*$-
algebra construction which is generated by $\X(L)$ in place of $\br_+$, but still modelled 
on Arveson's Wiener-Hopf approach. It appears that such an $\X(L)$-fractal-based 
$C^*$-algebra will serve as a spectral-invariant for our $B/L$ Hadamard-systems which 
are only symmetric, but generally {\it not\/} selfadjoint (relative to some choice of lattice 
$K$,  see Section 7 above).

The spectral-invariant question is an important one, and in our case we produce the dual 
representation pair (6.4) and (6.5) as a candidate. But representations $\{S_b\}$ of 
$\cO_N$ in the form (6.4), without a paired dual representation $\{T_\ell\}$, cf. (6.6), are 
present for iteration systems which are much more general than the affine fractals studied 
here. As a case in point we mention Matsumoto's \cite{Mat} recent analysis of (von 
Neumann type) cellular automata (details in \cite{MOW} and \cite{Od}); it is based on an 
$S$-representation which is given by a formula similar to our (6.4) above. There is also an 
associated endomorphism with an entropy that can be computed; but we stress that for 
these (and many other) iteration systems, there is typically {\it not\/} a dualitly based on 
exponentials $e^{i\lam \cdot x}$ and typically {\it not\/} a second $\{T_\ell \}$-
representation such that the two form a dual pair in any natural way.

We have studied the class of {\it spectral systems\/} $s=(R,B,L,K)$ in $\br^n$ with 
special view to the {\it selfadjoint\/} ones which are also of {\it Hadamard type\/}, see 
Lemma 4.1. (When $s$ is given in this class, the two sets $B$ and $L$ then have the same 
cardinality; it will be denoted $N$ for convenience in the following comments.) It is 
important (but elementary) {\it that this class of systems is closed under the tensor-
product operation\/}; i.e., if $s_1$ and $s_2$ are systems in $\br^{n_1}$ and $\br^{n_2}$ 
respectively, then the two properties (selfadjointness and Hadamard type) carry over to 
the system $s_1 \otimes s_2$ in $\br^{n_1 + n_2}=\br^{n_1} \times \br^{n_2}$. 

If the Hilbert spaces for the respective systems are $L^2(\mu_i)$, $i=1,2$; then the 
Hilbert space for $s_1 \otimes s_2$ is $L^2(\mu_1 \otimes \mu_2)$, and the measure 
$\mu_1 \otimes \mu_2$ is the unique probability measure on $\br^{n_1} \times 
\br^{n_2}$ which scales the affine tensor operations of $s_1 \otimes s_2$, see (2.3) 
above. The set $B$ for $s_1 \otimes s_2$ is $B_1 \times B_2$, and the matrix-operation 
is, $(b_1,b_2) \mapsto (R_1b_1,R_2b_2)$. In verifying the Hadamard property (3.2) for 
$s_1 \otimes s_2$, we use the important (known) fact that the class of generalized 
Hadamard matrices is closed under the tensor-product operation, i.e., that $U_1 \otimes 
U_2$ satisfies (3.2) with order $N=N_1 N_2$ if the individual factors $U_i$, $i=1,2$, do 
with respective orders $N_i$, $i=1,2$.

We say that a system $s$ is {\it irreducible\/} if it does not factor ``non-trivially'' $s\simeq 
s_1 \otimes s_2$; and we note that the examples above from Section 7 are all irreducible 
in this sense. (In fact this irreducibility notion is {\it different\/} from that of Section 6, but 
the examples are irreducible in both senses.)

The spectral geometry for regions in $\br^n$ has a long history, see e.g., 
\cite{Bo-Gu}, \cite{CV}, \cite{Ge}, and \cite{Gu-St}. But, so far, the {\it Laplace 
operator\/} has played a favored role despite the known incompleteness for the 
correspondence between the geometry of the given domain and the spectrum of the 
corresponding Laplace operator. The approach in \cite{De} is based instead on a 
multitude of {\it second order\/} differential operators, but the spectral correspondence is 
still incomplete there. Our present approach leads to a complete spectral picture and is 
based instead on a system of {\it first order\/} operators. For the fractal case however, the 
differential operators have no analogue.

While our simultaneous eigenfunctions are based, at the outset, on a commutative 
operator system, our spectral invariant derives instead from a dual pair of representations 
of a certain {\it non-abelian\/} (in fact {\it simple\/}) $C^*$-algebra.

Self-similar limit constructions have received much recent attention, starting with 
\cite{Hu}, and then more recently, see e.g., \cite{Ba-Gr}, \cite{Ed}, \cite{Ma}, and 
\cite{Ke}. These results seem to stress the geometry and the combinatorics of the infinite 
limits, and not the spectral theory. Our present emphasis is a {\it direct spectral/geometry-
correspondence\/}; and we also do {\it not\/} in \cite{JP6} impose the strict {\it 
expansivity\/} assumption (which has, so far, been standard almost everywhere in the 
literature). Furthermore, we wish to stress that the sets $\Om \subset \br^n$ which occur 
in our present spectral pairs {\it are more general\/} than the {\it self-reproducing tiles\/} 
(SRT) which were characterized in \cite{Ke\rm, Theorem 10}. However, Kenyon's SRT's 
can be shown to satisfy our conditions, although our class is {\it properly larger\/}; not 
only because of the expansivity assumption, but also because of the combinatorics, see 
\cite{Jo-Pe5} for details. Further work on these interconnections is also in progress.

\head Acknowledgments \endhead
Both authors were supported in part by grants from the U.S. National Science 
Foundation. It should finally be mentioned that the portion of the problem which related to 
commuting operator extensions for the partial derivatives ${\partial \over \partial x_j}$ 
originates with suggestions made first in 1958 by Professor I.E. Segal, see \cite{Fu} for 
details on this point. Encouragement and several correspondences from Professors B. 
Fuglede and R. Strichartz are also greatly appreciated. Detailed suggestions from 
Strichartz led to substantial improvements.

\enddocument